\documentclass[11pt]{article}
\usepackage{amssymb}
\usepackage{multicol}        
\usepackage[bottom]{footmisc}
\usepackage{amssymb}
\usepackage{amsmath}
\usepackage{mathrsfs}

\usepackage{caption}
\usepackage{times}
\usepackage[T1]{fontenc}
\usepackage{amsmath, amsthm, amsfonts}
\usepackage{amssymb}
\usepackage{amsbsy}
\usepackage{amscd}
\usepackage{verbatim}
\usepackage{xcolor}

\usepackage{epsfig}

\usepackage{placeins}

\title{\LARGE Estimates on the velocity of a rigid body moving in a fluid} 

\author{\Large Stathis Filippas,~~  Alkis Tersenov
 \\
                                                                           \\
        Department Mathematics and Applied Mathematics \\
         University of Crete,  \\
         70013 Heraklion,  Greece \\
        filippas@uoc.gr,~~~ tersenov@uoc.gr \\                                        
\\
    \\ }

\date{\today}

\begin{document}

\maketitle


\newcommand{\be}{\begin{equation}}
\newcommand{\ee}{\end{equation}}
\newcommand{\bea}{\begin{eqnarray}}
\newcommand{\eea}{\end{eqnarray}}
\newcommand{\bean}{\begin{eqnarray*}}
\newcommand{\eean}{\end{eqnarray*}}
\newcommand{\la}{\label}
\newcommand{\xa}{\alpha}
\newcommand{\xb}{\beta}
\newcommand{\xg}{\gamma}
\newcommand{\xG}{\Gamma}
\newcommand{\xd}{\delta}
\newcommand{\xD}{\Delta}
\newcommand{\xe}{\varepsilon}
\newcommand{\xz}{\zeta}
\newcommand{\xh}{\eta}
\newcommand{\vth}{\vartheta}
\newcommand{\Th}{\Theta}
\newcommand{\xk}{\kappa}
\newcommand{\xl}{\lambda}
\newcommand{\xL}{\Lambda}
\newcommand{\cF}{\mathcal{F}}
\newcommand{\cg}{\mathcal{G}}
\newcommand{\cR}{\mathcal{R}}
\newcommand{\cC}{\mathcal{C}}
\newcommand{\bn}{\boldsymbol{n}}
\newcommand{\bt}{\boldsymbol{\tau}}
\newcommand{\ber }{\boldsymbol{e_{r}}}
\newcommand{\bez }{\boldsymbol{e_{x_3}}}
\newcommand{\bet }{\boldsymbol{e_{\theta}}}
\newcommand{\bug }{\boldsymbol{e_{\xg}}}
\newcommand{\bxi}{\boldsymbol{\xi}}
\newcommand{\T}{\mathbb{T}}
\newcommand{\I}{\mathbb{I}}
\newcommand{\ba}{\boldsymbol{a}}
\newcommand{\bx}{\boldsymbol{x}}
\newcommand{\bu}{\boldsymbol{u}}
\newcommand{\bv}{\boldsymbol{v}}
\newcommand{\bg}{\boldsymbol{g}}
\newcommand{\bom}{\boldsymbol{\omega}}
\newcommand{\bF}{\boldsymbol{F}}
\newcommand{\bI}{\boldsymbol{I}}

\newcommand{\grad}{{\rm grad \,}}
\newcommand{\dive}{{\rm div \,}}

\newcommand{\xs}{\sigma}
\newcommand{\xS}{\Sigma}
\newcommand{\ps}{\psi}
\newcommand{\PS}{\Psi}
\newcommand{\xo}{\omega}
\newcommand{\xO}{\Omega}
\newcommand{\R}{ I \! \! R}
\newcommand{\ra}{\rightarrow}
\newcommand{\rft}{\rightarrow +\infty}
\newcommand{\finedim}{{\hfill $\Box$}}
\newcounter{newsection}
\newtheorem{theorem}{Theorem}[section]
\newtheorem{lemma}[theorem]{Lemma}
\newtheorem{prop}[theorem]{Proposition}
\newtheorem{coro}[theorem]{Corollary}
\newtheorem{defin}[theorem]{Definition}
\newcounter{newsec} \renewcommand{\theequation}{\thesection.\arabic{equation}}

\begin{abstract}

 We obtain estimates of all components of the velocity of a 3D rigid body moving in a viscous incompressible fluid  without any symmetry restriction on the shape of the rigid body or the container. The estimates are in terms of suitable norms of the velocity field in a small domain of the fluid only, provided the distance $h$ between the rigid body and the container is small.
  As a consequence we obtain suitable differential inequalities that control the distance $h$.
 The results are obtained 
 using  the fact that the vector field under investigation belongs to suitable
   function spaces, without any use of hydrodynamic equations.
 
\end{abstract}

\noindent {\bf AMS Subject Classification: } 46E35, 76D99   

\noindent {\bf Keywords: } Sobolev spaces, fluid solid interaction, non zero speed collision, non collision

\setcounter{equation}{0}
\section{Introduction}

Throughout the paper $\xO \subset \R^3$ is a domain,  $S \subset \xO$ is a bounded connected  domain
and $F = \xO \setminus\bar{S}$.   We denote by
 $ W^{1,p}_0(S, \,\xO)$, $p \geq 1$, the vector function space
consisting of functions
 \[
  \bu: \xO \to \R^3, ~~~~~~\bu \in W^{1,p}_0(\xO)  \ ,
 \]
 such that for a constant vector $\bu_* \in \R^3$ and a vector $\bom \in \R^3$ there holds 
\be\la{rig}
\bu(\bx)=\bu_{*}+ \bom \times  (\bx- \bx_{*}), ~~\hbox{ for}~~ \bx  \in S  \ ,
\ee
where  $\bx_{*} \in \overline{S}$ is a fixed point.

Such  function spaces arise naturally when studying the  motion of a rigid body $S$ inside a fluid region $\xO$,  see e.g., \cite{ HS,HT1, LT, V1}.
 In this context $\bu$ is the velocity field in $\xO$ where in particular inside the rigid body the velocity is given by (\ref{rig}).

We denote by $h$ the distance between $S$ and $\partial \xO$, that is,
\[
h =  {\rm dist}(S, \partial \xO) = |PQ| ,
\]
where  $P \in \partial S$ and $Q \in \partial \xO$ are two points that realize the distance. We are interested in the case where $h$ is small, hence we may assume that
 \[
h < H:= {\rm diam} S .
\] 
We introduce an orthogonal coordinate system $(x_1, x_2, x_3)$ with  the origin at the point $Q$. We note that
the plane $\bx' = (x_1, x_2)$ is tangent to $\partial \xO$ at the point $Q$ and parallel to the tangent plane of $\partial S$ at the point $P=(0,0,h)$.
Let $\bu_p = (\bu_{P \tau}, u_{P3})=(u_{P1}, u_{P2},  u_{P3})$ be the velocity of  the rigid body at the point 
$P$ and $\bom= (\bom_{\tau}, \xo_3)= ( \xo_1, \xo_2, \xo_3)$ its angular velocity.


 \setlength{\unitlength}{1,4cm}
 
  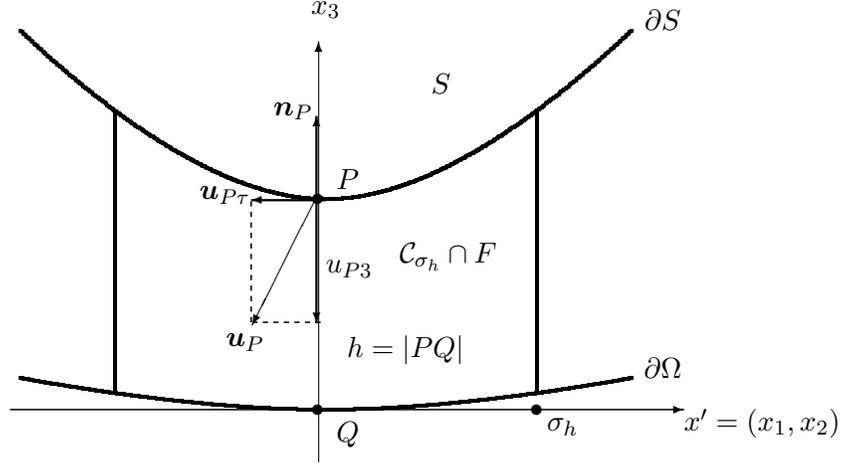
\begin{figure}
 
 \begin{picture}(5,4)(-5.2, -2.6)

\put(0,-2.5){\vector(0,1){4}}.                 
\put(0,1.8){\makebox(0,0){$x_3$}}

\linethickness{0.4mm}
\qbezier(0.0,0.0)(1.2,0.0).  
(2.9,1.6)
\linethickness{0.4mm}
\put(3.2 ,1.7){\makebox(0,0){$\partial S$}}

\linethickness{0.4mm}
\qbezier(0.0,0.0)(-1.2,0.0).  
(-2.9,1.6)

\linethickness{0.4mm}
\qbezier(0,-2)(1.2,-2).  
(2.9, -1.7)

\put(3.2 ,-1.6){\makebox(0,0){$\partial \Omega$}}

\linethickness{0.4mm}
\qbezier(0, -2)(-1.2,-2).  
(-2.9, -1.7)

\put(-0.08,0){\circle*{0.1}}
\put(-0.08,-2){\circle*{0.1}}
\put(2,-2){\circle*{0.1}}


\put(-0.1,0){\vector(-1,-2){0.6}}

\linethickness{.2mm}
\multiput(-0.1,-1.17)(-.1,0){7}{\line(1,0){.05}}

\multiput(-0.71,-0.07)(0,-0.1){12}{\line(0,1){.05}}

\linethickness{.3mm}
\put(-0.,-0.01){\vector(-1,0){.71}}
\put(-0.09,-0.01){\vector(0,-1){1.17}}

\put(-0.95, -1.4){$\bu_{P }$}
\put(-1.2, 0){$\bu_{P \tau }$}
\put(-0., -0.7){$u_{P 3}$}

\put(0.7, -0.6){$\cC_{\xs_h} \cap F$}


\linethickness{.3mm}
\put(-0.09,-0.01){\vector(0,1){0.8}}
\put(-.5,0.8){$\bn_P$}

\linethickness{.4mm}
\put(-2,-1.82){\line(0,1){2.67}}

\put(2,-1.82){\line(0,1){2.67}}

\linethickness{0.2mm}
\put(-3,-2){\vector(1,0){6.4}}     
\put(3.4, -2,2){$x'=(x_1,x_2)$}

\put(1.,  1.){$S$}
\put(0.1,  0.1){$P$}
\put(0.1,  -2.3){$Q$}
\put(0.2, -1.5){$h=|PQ|$}
\put(2.1, -2.2){$\sigma_h$}

\end{picture}

\caption{ The coordinate system}

 \end{figure}

Our first result concerns vector fields that enjoy the regularity of weak solutions of the solid--fluid interaction problem in hydrodynamics cf \cite{HT1, HS, SST}.
 For 
\[
\cC_{\xs} \cap F = \{|\bx'|<\xs \} \cap F,
\]
 we have
\begin{theorem} \la{wk} 
We  assume that 
both $\xO$ and $S$ are uniformly $C^{1,\xa}$ domains, $\xa \in (0, 1]$.
Let $\bu \in W^{1,p}_0(S, \,\xO) $ be such that ${\rm div} \bu = 0$. Then, there exist positive constants $H_{*}$, $C_w$, $c_0$ depending only on $S$, $\xO$, $p$, such that, whenever $h < H_{*}$ and $\xs_{h}:= c_0 h^{\frac{1}{1+\xa}}$,
\bean
  |u_{P3}| & \leq &  \, C_w  h^{\frac{1+2 \xa}{p(1+\xa)}
  \left( p - \frac{3+\xa}{1+ 2 \xa} \right)}  
\| \nabla \bu \|_{L^p(\cC_{\xs_h}\cap F)},  \\
  |\bu_{P\tau}| & \leq  & C_w  \, h^{\frac{1}{p}
  \left( p - \frac{3+\xa}{1+ \xa} \right)}  \| \nabla \bu \|_{L^p(\cC_{\xs_h}\cap F)},  \\
  |\bom_{\tau}|  & \leq  &  C_w   \, h^{\frac{2 \xa}{p (1+\xa)}
  \left( p - \frac{3+\xa}{2 \xa} \right)}  \| \nabla \bu \|_{L^p(\cC_{\xs_h}\cap F)} , \\
  |\xo_3|  & \leq  &  C_w  \,  h^{\frac{ \xa}{p (1+\xa)}
  \left( p - \frac{3+\xa}{ \xa} \right)}  \| \nabla \bu \|_{L^p(\cC_{\xs_h}\cap F)} .
  \eean
\end{theorem}
Estimates for $|u_{P3}|$ where first  proved in  \cite{HS, V1}  where in the right hand side the norm
of $\nabla \bu $ over a larger region is required. The novelty then in our Theorem is that the norm 
of $\nabla \bu $  appears in a small region of the {\em fluid only}. This is in agreement with the observation that,
in simple geometries,  knowledge of the  fluid vector field in 
a small region of the fluid between $S$ and $\partial \xO$ is enough to produce significant information such  as the noncollision of the rigid body and the boundary of the container, see e.g. \cite{GH1, H, HT1, LG}. The rest of the estimates are altogether new.  Since the velocity at any point of $S$ is given by (\ref{rig}), one obtains a complete  control of the motion of the rigid body as it approaches the boundary of $\xO$.
The above result  extends  Theorem 4.6  of \cite{FT} in 3 dimensions.

Our next result concerns vector fields that enjoy the regularity of strong solutions of the solid--fluid interaction problem, cf  \cite{Tak,HT2}.

\begin{theorem} \la{str} 
We  assume that 
both $\xO$ and $S$ are uniformly $C^{3}$.
Let $\bu \in W^{1,2}_0(S, \,\xO) $ be such that ${\rm div} \bu = 0$ and in addition $ |\nabla^2 \bu |\in L^{2}(F)$ .
Then, there exist positive constants $H_{*}$, $C_s$, $c_0$ depending only on $S$, $\xO$, such that, whenever $h < H_{*}$ and $\xs_{h}:= c_0 h^{\frac{1}{2}}$,
\[
|u_{P3}| \leq   C_{s} \left( h^{\frac32} \,
 \| \nabla^2 \bu\|_{L^2( \mathcal{C}_{\xs_h}\cap F)} +  h \,
 \| \nabla \bu\|_{L^2( \mathcal{C}_{\xs_h}\cap F)}  \right) .
 \]

\end{theorem}
This can be thought as an improvement of the corresponding estimate of Theorem \ref{wk}, which for $\xa=1$, $p=2$ takes the form
\[
|u_{P3}| \leq   C_{w} h^{\frac12} \,
 \| \nabla \bu\|_{L^2( \mathcal{C}_{\xs_h}\cap F)} .
 \]
Such estimates for strong solutions were first established in \cite{H05, HT2} in the case where $\xO$ is half--space and $S$ is a disc or a ball. In such a case only the first term in the right hand side appears due to the apparent symmetries. In the absence of symmetries however, the second term above seems to be necessary. The proof of Theorem \ref{str} uses in an essential way all the estimates of Theorem \ref{wk}.
It is not clear whether  the other quantities appearing in Theorem \ref{wk} can also be improved in a similar way
and we leave it as an open question.

We next consider the dynamic case, where $S$ is allowed to move inside $\xO$.
 Denote by $S(t) \subset \xO$ the position of $S$ at time $t$ and $F=F(t)= \xO \setminus S(t)$. 
  We thus assume that there are $L^{\infty}(0,T)$ functions $\bu_{*}(t)$ and $\bom(t)$ so that 
$S(t)$ consists of points $\bx(t)$  satisfying
\[
\frac{d \bx(t)}{ d t} = \bu_{*}(t) + \bom(t) \times (\bx(t) - \bx_{*}(t)), ~~~~~~\bx(0)= \bx_0 \in S_{0}.
\]
For $p \geq 1$, $q \geq 1$ we define the   space of functions $\bu (\bx, t)$, $\bx \in \xO$, $t \in (0,T)$,
\bean
L^q(0,T; W^{1,p}_0(S(t) \, , \xO)) :=  &  &
 \\
  &  &  \hspace{-5,3cm}
\{ \bu \in   L^q \left(0,T;  W^{1,p}_0(\xO) \right):
~
  \bu(\bx,t) 
   =\bu_{*}(t)+ \bom(t) \times  (\bx- \bx_{*}(t)) ~\mbox{ for}~ \bx  \in S(t)  
 \} .
 \eean
The distance between $S(t)$ and $\partial \xO$ is now  a function of time $h(t)$,  that is,
\[
h(t) =  {\rm dist}(S(t), \partial \xO) .
\]
The following is a  consequence of the previous Theorems.

\begin{theorem} \la{dynam} 
\noindent ({\rm i}) We assume that both $\xO$ and $S_0$ are uniformly $C^{1,\xa}$, $\xa \in (0,1]$ domains.
In addition $\xO$ satisfies an inner sphere condition. Let
\be\la{ws}
\bu \in L^{\infty}(0,T; L^{2}(\xO)) \cap L^{1}(0,T; W^{1,p}_0(S(t),\xO)),~~ p \in [1, \infty),~~ {\rm div} \bu = 0.
\ee
Then, there exist positive constants $H_{0}$, $C_w$, $c_0$ depending only on $S$, $\xO$, $p$, such that, whenever $h < H_{0}$, function $h=h(t)$ is Lipschitz continuous and for $\xs_{h}:= c_0 h^{\frac{1}{1+\xa}}$ there holds
\be\la{th19}
\left| \frac{d h}{dt} \right| \leq C_w \,  h^{\frac{1+2 \xa}{p(1+\xa)}
  \left( p - \frac{3+\xa}{1+ 2 \xa} \right)}  
\| \nabla \bu \|_{L^p(\cC_{\xs_h}\cap F)}, ~~~~\mbox{for a.a.}~~ t \in (0,T).
 \ee
 
\noindent ({\rm ii}) We assume that both $\xO$ and $S_0$ are uniformly $C^{3}$  domains. The vector 
field $\bu$  satisfies (\ref{ws}) with $p=2$ and in addition
\[
| \nabla^2 \bu |\in   L^{2}(F(t)), ~~~~~ \mbox{ for a.a.}~t \in (0,T).
\]
Then, there exist positive constants $H_{0}$, $C_s$, $c_0$ depending only on $S$, $\xO$,  such that, whenever $h < H_{0}$, function $h=h(t)$ is Lipschitz continuous and for $\xs_{h}:= c_0 h^{\frac{1}{2}}$,
 \be\la{s4}
\left| \frac{d h}{dt} \right| \leq C_{s} \left( h^{\frac32} \,
 \| \nabla^2 \bu\|_{L^2( \mathcal{C}_{\xs_h}\cap F)} +  h \,
 \| \nabla \bu\|_{L^2( \mathcal{C}_{\xs_h}\cap F)}  \right) ,~~~~\mbox{for a.a.}~~ t \in (0,T).
 \ee
\end{theorem}
As a consequence of Theorem \ref{dynam}(ii),  whenever $\lim_{t \to T}h(t) \to 0$, one has that
\[
\lim_{t \to T}\int_{0}^{t}  h^{\frac12}(\tau) \|\nabla^2 \bu \|_{L^2( \mathcal{C}_{\xs_h(\tau)}\cap F(\tau))} \, d \tau = \infty ,
\]
 see Corollary \ref{cor}.
 
 The focus of the paper is in the 3D case, but similar results hold for the 2D case as well. We present some of them in the form of Remarks.
  We also note that although  our results are obtained for a solenoidal vector  field, a similar approach can be applied in the nonsolenoidal   case.
 
The paper is organized as follows. In section 2 the proof of Theorem \ref{wk} is given and in section 3 the proof of Theorems \ref{str} and \ref{dynam} are given. In section 4 we present a 3D  example showing the optimality of Theorem \ref{wk}. In addition, in this example a smooth ($\xa=1$) rigid body hits the boundary of the container with nonzero speed.

\setcounter{equation}{0}
\section{Estimates for weak solutions }
Here we will give the proof of Theorem \ref{wk}. 
 We first recall   that a domain $D \subset \R^3$ is {\em uniformly $C^{1,\xa}$}, $\xa \in (0,1]$,  if there exist constants $k>0$ and $\xs_0>0$ s.t. for each $x_0 \in \partial D$ there is an isometry $\cR: \R^3 \to \R^3$ s.t.
 \[
 \cR(0) = x_0, ~~~~\cR(\{x\in \R^3: 
 ~k|\bx'|^{1+\xa} <|x_3| < k|\xs_0|^{1+\xa} , ~~~|\bx'|<\xs_0 \}) \cap  \partial D = \emptyset.
 \]
We note that if $D$ is a bounded  $C^{1,\xa}$ domain then it is a 
  uniformly $C^{1,\xa}$ domain, see e.g, \cite{AB}.

   We will first prove a variation of Theorem \ref{wk}. To this end
we fix a  $\xs \in [0, \xs_0/2]$ and let $0<   \rho \leq \xs$. We denote by $\mathcal{C}_{\rho}$ the cylinder
 of radius $\rho$ and height $h+3k \xs^{1+\xa}$,
 \bean
 \mathcal{C}_{\rho} & = & \{(r, \theta, x_3):~  0 < r < \rho,  ~~0 \leq \theta < 2 \pi, ~-k \xs^{1+\xa} < x_3 < h+2k \xs^{1+\xa}
  \}.  
 \eean

 \begin{theorem} \la{th1}
 We  assume that 
both $\xO$ and $S$ are uniformly $C^{1,\xa}$ domains, $\xa  \in (0,1]$.
Let $\bu \in W^{1,p}_0(S, \,\xO) $ be such that ${\rm div} \bu = 0$. Then, there exist positive constant $C_w$ depending only on $S$ and $\xO$, such that if
 $h \in [0,H]$  and  $\xs \in (0, \xs_0/2]$   there holds
\bean
| u_{P3} | & \leq  & C_w \xs^{-1-\frac{2}{p}} (h + 3 k \xs^{1+\xa})^{2-\frac{1}{p}}  
\| \nabla \bu \|_{L^p(\cC_{\xs})},   \\
 |\bu_{P\tau}| & \leq  & C_w  \xs^{-1-\xa-\frac{2}{p}} (h + 3 k \xs^{1+\xa})^{2-\frac{1}{p}}
\| \nabla \bu \|_{L^p(\cC_{\xs})} , \\
|\bom_{\tau}| &  \leq  & C_w  \xs^{-2-\frac{2}{p}} (h + 3 k \xs^{1+\xa})^{2-\frac{1}{p}}
\| \nabla \bu \|_{L^p(\cC_{\xs})}, \\
 |\xo_{3}| &  \leq  & C_w  \xs^{-2-\xa-\frac{2}{p}} (h + 3 k \xs^{1+\xa})^{2-\frac{1}{p}}
\| \nabla \bu \|_{L^p(\cC_{\xs})}.
\eean

%
%
%
%
 \end{theorem}

 \noindent {\em Proof:}  Let $0<   \rho \leq \xs  \leq \xs_0/2$. We note that $\xs$ will be fixed throughout the proof. In addition we assume
 that $0< \xg  < \xd  \leq \pi$.  We first fix some notation.

 \vspace{3mm} \noindent  \underline{\bf Notation.} 
 We denote by $D_{\rho}$ the disc of radius $\rho$ and by $D_{\rho, \xg}$ the half disc given by
 \[
 D_{\rho, \xg} := \{  0 < r < \rho,  ~~\xg \leq \theta < \pi+ \xg \}.
 \]
 We also denote by  $\mathcal{C}_{\rho, \xg } $ the half  cylinder given by
 \bean
 \mathcal{C}_{\rho, \xg }  =  \{(r, \theta, x_3):~  (r,~ \theta) \in  D_{\rho, \xg}, ~-k \xs^{1+\xa} < x_3 < h+2k \xs^{1+\xa}  \}.  
 \eean
Clearly we have $\mathcal{C}_{\rho} = \mathcal{C}_{\rho, \xg } \cup \mathcal{C}_{\rho, \pi +\xg} $.
We also set \[
\xd \mathcal{C}_{\rho, \xg } := \mathcal{C}_{\rho, \xg }
\setminus \mathcal{C}_{\rho/2, \xg } .
\]  
 To describe  the boundary  surfaces of $\mathcal{C}_{\rho, \xg }$ we use the notation
(upper, lateral curved, lateral flat,  lower)
 \bean
 \Gamma_{\rho, \xg}^{+} & = & \{(r, \theta, x_3):~(r,~ \theta) \in  D_{\rho, \xg},  ~~  x_3 = h+2k \xs^{1+\xa} \} , \\
 \Gamma_{\rho,\xg}^{0} & = & \{(r, \theta, x_3):~   r = \rho,  ~~ \xg < \theta <  \pi+ \xg, ~~ -k \xs^{1+\xa} <  x_3 < h+2k \xs^{1+\xa} \} , \\
 \Gamma_{\rho,\xg}^{1} & = & \{(r, \theta, x_3):~   0<r < \rho,  ~ ~\theta = \xg, ~~ -k \xs^{1+\xa} <  x_3 < h+2k \xs^{1+\xa} \} , \\
 \Gamma_{\rho, \xg}^{-} & = & \{(r, \theta, x_3):~(r,~ \theta) \in  D_{\rho, \xg},  ~~  x_3 = - k \xs^{1+\xa} \}  .
 \eean
 The (total) lateral  surface of $\mathcal{C}_{\rho, \xg }$ is given by
 \[
 \Gamma_{\rho, \xg} := \Gamma^{0}_{\rho, \xg} \cup \Gamma^{1}_{\rho, \xg} \cup \Gamma^{1}_{\rho, \pi+\xg} .
 \]
 Concerning $\xd \mathcal{C}_{\rho, \xg } $ we have that the upper surface is given by
 \[
 \xd \Gamma_{\rho, \xg}^{+} := \Gamma_{\rho, \xg}^{+} \setminus
 \Gamma_{\rho/2, \xg}^{+}.
 \]
 The (total) lateral  surface of $\xd \mathcal{C}_{\rho, \xg } := \mathcal{C}_{\rho, \xg }
\setminus \mathcal{C}_{\rho/2, \xg }$ is given by
\[
\xd \Gamma_{\rho, \xg} := \Gamma^{0}_{\rho, \xg} \cup \Gamma^{0}_{\rho/2, \xg}  \cup \xd \Gamma^{1}_{\rho, \xg} \cup \xd \Gamma^{1}_{\rho, \pi+\xg},
\]
where
\[
\xd \Gamma^{1}_{\rho, \xg}  := \Gamma^{1}_{\rho, \xg} \setminus \Gamma^{1}_{\rho/2, \xg}  .
\]

\begin{figure}
\begin{center}
\includegraphics[scale=0.7]{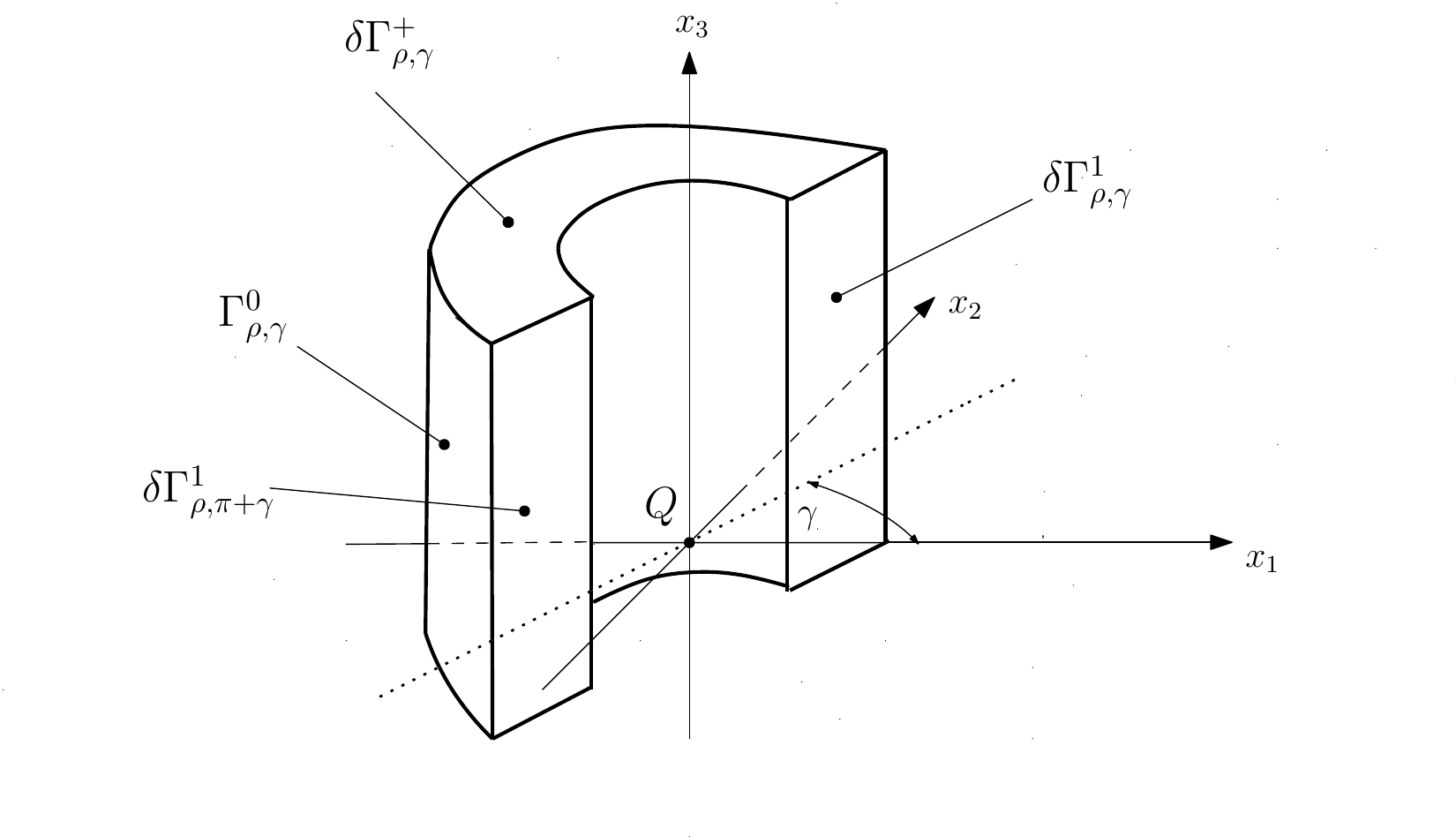}
\caption{ The domain $\xd \mathcal{C}_{\rho, \xg }$ 
 and its boundary surfaces.}
\end{center}
\end{figure}

\vspace{3mm} \noindent \underline{\bf Plan of the proof.}
The estimate of $|u_{P3}|=| \bu_P \cdot \bn_P |$ is essentially  given in \cite{V1} Theorem 3.1, see pp 321-322.
 In the sequel we will first estimate the tangential part of 
$\bom$, that is, $\bom_{\bt}=(\xo_1, \xo_2,0)$, then the tangential part of 
$\bu_P $, that is $ \bu_{P \tau}=(u_{P1},u_{P2},0)$ and finally the third component of $\bom$, that is 
 $|\xo_3|$. We follow this order in the proof,  because at each step we use the  estimates of the previous steps.

 \vspace{3mm} \noindent \underline{\bf Estimate of $\bom_{\bt}$.} We will apply 
 divergence Theorem in the domain $\xd \mathcal{C}_{\rho, \xg } = \mathcal{C}_{\rho, \xg } \setminus   \mathcal{C}_{\rho/2, \xg }$. By choosing this non symmetric domain the 
 tangential component of $\bom$ appears in the calculations. The removal of $\mathcal{C}_{\rho/2, \xg }$  is done for technical reasons and will be explained later in the proof.

Since $\bu$ is divergence free  and $\bu=0$ on $ \Gamma_{\rho, \xg}^{-} $, applying the divergence Theorem in  $\xd \mathcal{C}_{\rho, \xg }$ we have that
\bea\la{d2om}
\int_{\xd \xG_{\rho, \xg }^{+}} \bu \cdot \bn \, dS + \int_{\xd \Gamma_{\rho, \xg} } \bu \cdot \bn \, dS   = 0.
 \eea
 We next compute the integral over $\xd \Gamma_{\rho, \xg}^{+} = \Gamma_{\rho, \xg}^{+} \setminus
 \Gamma_{\rho/2, \xg}^{+}$, where we have
  that $\bu = \bu_P +  \bom \times (\bx - \bx_P)$  and $ \bn= \bn_p =(0,0,1)$. A
 straightforward calculation yields
 \bean
 \int_{ \xd \Gamma_{\rho, \xg}^{+}} \bu \cdot \bn \, dS &  =  &
 \int_{ \xd \Gamma_{\rho, \xg}^{+}} \bu_P \cdot \bn_P \, dS +  \int_{\xd \Gamma_{\rho, \xg}^{+}}  \bom  \times
 (\bx - \bx_P) \cdot  \bn \, dS  \\  
  &  =  & \frac38 \pi \rho^2 \bu_P \cdot \bn_P +
\bom \cdot   \int_{ \xd \Gamma_{\rho, \xg}^{+}}   (\bx - \bx_P) \times  \bn \, dS  \\
 &  =  & \frac38 \pi \rho^2 \bu_P \cdot \bn_P +  \frac{7 \rho^3 }{12} (\xo_1 \, \cos \xg +  
\xo_2 \, \sin \xg  ).
 \eean
 From this and (\ref{d2om}) we get
 \be\la{a3}
 \frac{7 \rho^3}{12}  |\xo_1 \,\cos \xg+ ~\xo_2 \, \sin \xg  | \leq \frac{3 \pi\rho^2 }{8}  | \bu_P \cdot \bn_P | +  \int_{ \Gamma_{\rho, \xg}^{0}\cup \Gamma_{\rho/2, \xg}^{0}} |\bu| dS 
 +\int_{ \delta \Gamma^{1}_{\rho, \xg} \cup \delta \Gamma^{1}_{\rho, \pi+\xg}} |\bu|\, dS .
 \ee
 To estimate the terms in the right hand side we  will integrate  the above estimate from  $\xg=0$ to $\xd \leq \pi$
 and then from
 $\rho=\xs/2$ 
  to $\xs$.  To  this end we
  first note that 
 \bean
 \int_{0}^{\xd} \int_{ \Gamma_{\rho, \xg}^{0} } |\bu|\,  dS  d \xg  & =  & \int_{0}^{\xd} \int_{\xg}^{\pi + \xg} \int_{-k \xs^{1+\xa}}^{h+2k \xs^{1+\xa}}         |\bu (\rho, \theta, x_3)| \,\rho  ~dx_3 d \theta d \xg  \\
 & \leq & \int_{0}^{\xd} \int_{0}^{2\pi } \int_{-k \xs^{1+\xa}}^{h+2k \xs^{1+\xa}}         |\bu (\rho, \theta, x_3)| \,\rho  ~dx_3 d \theta d \xg  \\
  & \leq &  \pi \int_{0}^{2\pi } \int_{-k \xs^{1+\xa}}^{h+2k \xs^{1+\xa}}         |\bu (\rho, \theta, x_3)| \,\rho  ~dx_3 d \theta.
 \eean
We next inegrate  from
 $\rho=\xs/2$ 
  to $\xs$ to get
  \bean
\int_{\xs/2}^{\xs} \int_{0}^{\xd} \int_{ \Gamma_{\rho, \xg}^{0} } |\bu|\,  dS  d \xg d \rho
 & \leq  & \pi \int_{\cC_{\xs} \setminus \cC_{\xs/2}} |\bu| ~d \bx \leq \pi \int_{\cC_{\xs} } |\bu| ~d \bx  \\
 & \leq  & \pi ~\| \bu \|_{L^p(\cC_{\xs})} ~|\cC_{\xs}|^{\frac{p-1}{p}}  \\
 & \leq  & \pi~ C_1 ~(h + 3 k \xs^{1+\xa}) ~ \|\nabla \bu \|_{L^p(\cC_{\xs})}~ (\xs^2 (h + 3 k \xs^{1+\xa}))^{\frac{p-1}{p}}  \\
 & =  &  \pi ~C_1 ~\xs^{2-\frac{2}{p}} ~(h + 3 k \xs^{1+\xa})^{2-\frac{1}{p}}~
\| \nabla \bu \|_{L^p(\cC_{\xs})},
 \eean
 for some universal positive constant $C_1$.
 Working similarly for the integral over $\Gamma_{\rho/2, \xg}^{0}$ we end up with
   \bean
\int_{\xs/2}^{\xs} \int_{0}^{\xd} \int_{ \Gamma_{\rho, \xg}^{0} \cup \Gamma_{\rho/2, \xg}^{0} } |\bu|\,  dS  d \xg d \rho
 & \leq  & 3 \pi ~\int_{\cC_{\xs} } |\bu| ~d \bx \\
  & \leq  &  3 \pi~ C_1 ~\xs^{2-\frac{2}{p}} ~(h + 3 k \xs^{1+\xa})^{2-\frac{1}{p}}~
\| \nabla \bu \|_{L^p(\cC_{\xs})}.
 \eean
 Concerning  the flat parts we have
 \bean
 \int_{0}^{\xd}  \int_{\xs/2}^{\xs}  \int_{ \xd \Gamma_{\rho, \xg}^{1}} |\bu|\, dS  d \rho d \xg &  = & \int_{0}^{\xd} 
 \int_{-k \xs^{1+\xa}}^{h+2k \xs^{1+\xa}}     \int_{\xs/2}^{\xs}        \int_{\rho/2}^{\rho}  |\bu(r, \theta, x_3)|  \, dr d \rho  dx_3 d \xg  \\
 & \leq &  \int_{0}^{\xd}   \int_{-k \xs^{1+\xa}}^{h+2k \xs^{1+\xa}}     \int_{\xs/2}^{\xs}        \int_{\xs/4}^{\xs}  |\bu|  \, dr d \rho  dx_3 d \xg  \\
  &  = & \frac{\xs}{2}   \int_{0}^{\xd}   \int_{-k \xs^{1+\xa}}^{h+2k \xs^{1+\xa}}            \int_{\xs/4}^{\xs} \frac{|\bu|}{r}   ~r \, dr   dx_3 d \xg  \\
   &  \leq  &  \frac{\xs}{2}
    \int_{\cC_{\xs} \setminus \cC_{\xs/4} }\frac{|\bu|}{r} d \bx \\
     &  \leq & 2 \int_{\cC_{\xs} } |\bu| ~d \bx \\
  & \leq  &  2 C_1 ~\xs^{2-\frac{2}{p}} ~(h + 3 k \xs^{1+\xa})^{2-\frac{1}{p}}~
\| \nabla \bu \|_{L^p(\cC_{\xs})}.
 \eean
 It is in the above calculation that the removal of of $\mathcal{C}_{\rho/2, \xg }$  is required.
The same estimate holds true for the integral over $\xd \Gamma_{\rho, \pi+\xg}^{1}$.

On the other hand examination of the proof of Theorem 3.1 in \cite{V1} (see p. 321-322) shows that the following estimate is true
\be\la{vic1}
|u_{P3}|=| \bu_P \cdot \bn_P | \leq C_2 \xs^{-1-\frac{2}{p}} (h + 3 k \xs^{1+\xa})^{2-\frac{1}{p}}
\| \nabla \bu \|_{L^p(\cC_{\xs})},
\ee
for a universal positive constant $C_2$.
    
 Integrating twice (\ref{a3}) and employing the above estimates we get
 \[
 |(\xo_1 \, \sin \xd  + \xo_2 \,(1 -\cos \xd)| \leq C_0  \xs^{-2-\frac{2}{p}} (h + 3 k \xs^{1+\xa})^{2-\frac{1}{p}}
\| \nabla \bu \|_{L^p(\cC_{\xs})}. 
\]
 Since $\xd \in (0, \pi)$ is arbitrary, we conclude
 \be\la{bf1}
|\bom_{\tau}| \leq  |\xo_1| + |\xo_2| \leq C  \xs^{-2-\frac{2}{p}} (h + 3 k \xs^{1+\xa})^{2-\frac{1}{p}}
\| \nabla \bu \|_{L^p(\cC_{\xs})} ,
\ee
 for a universal positive constant $C$.

\vspace{5mm}
   
    \noindent \underline{{\bf Estimate of $\bu_{P\tau}$.}} We will use a similar argument, but this time instead of working with the half cylinder  $\mathcal{C}_{\rho, \xg } $  we will work with a half cylindrical domain with $\cup$--like top  surface given by

\[
\mathcal{C}_{\rho, \xg }^{\cup} := \{(r, \theta, x_3):~  (r,~ \theta) \in  D_{\rho, \xg}, ~-k \xs^{1+\xa} < x_3 < h+ k r^{1+\xa}  \}.  
 \]
Such a choice allows for the tangential part of     $\bu_{P}$ to enter the calculations.

To describe  the boundary  surfaces of $\mathcal{C}_{\rho, \xg }^{\cup}$ we use, in complete analogy with the previous case,  the notation
 \[
 \Gamma_{\rho, \xg}^{\cup,+}  =  \{(r, \theta, x_3):~(r,~ \theta) \in  D_{\rho, \xg},  ~~  x_3 = h+ k r^{1+\xa} \} , 
 \]
 for the upper surface and quite similarly $\Gamma_{\rho,\xg}^{\cup, 0} $, 
 $\Gamma_{\rho,\xg}^{\cup, 1}$, $ \Gamma_{\rho, \xg}^{\cup, -}=   \Gamma_{\rho, \xg}^{ -}  $ for the 
  lateral curved, lateral flat,  and lower surfaces respectively.
 Again, we apply divergence  Theorem in $ \xd \mathcal{C}_{\rho, \xg }^{\cup}:=\mathcal{C}_{\rho, \xg }^{\cup} \setminus   \mathcal{C}_{\rho/2, \xg }^{\cup} $ to get the analogue of (\ref{d2om})
 \bea\la{d3om}
\int_{ \Gamma_{\rho, \xg}^{\cup,+}\setminus \Gamma_{\rho/2, \xg}^{\cup,+} } \bu \cdot \bn \, dS + \int_{ \Gamma_{\rho, \xg}^{\cup,0} \cup \Gamma_{\rho/2, \xg}^{\cup,0}   } \bu \cdot \bn \, dS  
 +\int_{\delta \Gamma^{\cup,1}_{\rho, \xg} \cup \delta \Gamma^{\cup,1}_{\rho, \pi+\xg}  } \bu \cdot \bn \, dS = 0.
 \eea
The first integral above is computed explicitly. In fact, straightforward
calculations show that 
\bea\la{cup+}
 \int_{ \Gamma_{\rho, \xg}^{\cup,+}} \bu \cdot \bn \, dS  & = &
  \frac12 \pi \rho^2 \bu_P \cdot \bn_P + \frac{2k(1+\xa)}{2+\xa} \rho^{2+ \xa} \bu_P \cdot (\sin \xg, -\cos \xg, 0) \nonumber  \\
  & & + \left( \frac23  +\frac{2 k^2 (\xa+1)}{2 \xa+3} 
  \rho^{2 \xa} \right)\,\rho^3 \, \bom \cdot (\cos \xg, \sin \xg, 0).
 \eea
For the remaining two integrals in (\ref{d3om}) we note that 
\[
\Gamma_{\rho,\xg}^{\cup, 0}
\subset \Gamma_{\rho,\xg}^{ 0}~~~~\mbox{and~ similarly }~~~~  \delta \Gamma^{\cup,1}_{\rho, \xg} \subset \delta \Gamma^{1}_{\rho, \xg}.
\]
From (\ref{d3om}) then, we obtain the analogue of (\ref{a3}) which reads
\bea\la{a32}
 A_1 \rho^{2+\xa} |u_{P1} \,\sin \xg - ~u_{P2} \, \cos \xg  | & \leq  &
  A_2 \rho^2 | \bu_P \cdot \bn_P | +  A_3  \rho^3 |\bom_{\tau}| \nonumber \\
  & & \hspace{-2cm}+ \int_{ \Gamma_{\rho, \xg}^{0}\cup \Gamma_{\rho/2, \xg}^{0}} |\bu| dS 
 +\int_{ \delta \Gamma^{1}_{\rho, \xg} \cup \delta \Gamma^{1}_{\rho, \pi+\xg}} |\bu|\, dS,
 \eea
 with positive constants $A_i$, $i=1,2,3$, depending only on $k$, $\xa$.
 Again,
  we  will integrate  the above estimate from  $\xg=0$ to $\xd \leq \pi$
 and then from
 $\rho=\xs/2$ 
  to $\xs$ to reach the analogue of (\ref{bf1}),
\be\la{bf2}
|\bu_{P\tau}| \leq  |u_{P1}| + |u_{P2}| \leq C  \xs^{-1-\xa-\frac{2}{p}} (h + 3 k \xs^{1+\xa})^{2-\frac{1}{p}}
\| \nabla \bu \|_{L^p(\cC_{\xs})} ,
\ee
with a positive constant $C=C(k,\xa)$.

 \begin{figure}
\centering
 
 \setlength{\unitlength}{1,2cm}

 \begin{picture}(12 ,3)(-5.5, -0)
 
\put(0.3,-4.5){\vector(0,1){7.4}}.                 
\put(0.2,3.2){\makebox(0,0){$x_3$}}

\linethickness{0.1mm}. 
\qbezier(0.,0)(1.5,0)            
(3.2,1)
\qbezier(0.0,0)(-1.5,0)
(-2.8,1)

\linethickness{0.1mm}.          
\qbezier(0.0, -3)(1.5, -3)
(3.1, -4)
\qbezier(0.0, -3)(-1.5, -3)
(-2.9,-4)

\color{green}                           
\linethickness{0.4mm}.  
\qbezier(0, 0.02)(1., -0)
(2, 0.49)
\qbezier(0, 0.02)(-1, -0)
(-2, 0.41)

\color{black}

\linethickness{0.4mm}
\qbezier(0,0)(1.2,0).  
(4,0.5)
\linethickness{0.2mm}
\put(4 ,0.7){\makebox(0,0){$\partial S$}}

\linethickness{0.4mm}
\qbezier(0.0,0.0)(-1.2,0.0).  
(-4,0.5)

\linethickness{0.4mm}
\qbezier(0,-3)(1.2,-3).  
(3.8, -2.7)

\put(3.8 ,-2.5){\makebox(0,0){$\partial \Omega$}}

\linethickness{0.4mm}
\qbezier(0, -3)(-1.2,-3).  
(-3.8, -2.7)

\put(-0.05,0){\circle*{0.1}}               
\put(-0.05,-3){\circle*{0.1}}
\put(-0.05,-4){\circle*{0.1}}
\put(2,-3){\circle*{0.1}}
\put(3,-3){\circle*{0.1}}
\put(-0.05,1){\circle*{0.1}}   
\put(-0.05,2){\circle*{0.1}}

\color{blue}
\linethickness{0.4mm}
\put(-2,-4){\line(0,1){6}}                
\put(2,-4){\line(0,1){6}}
\put(-2,-4){\line(1,0){4}}
\put(-2,2){\line(1,0){4}}

\color{red}

\put(-2,1.153){\line(6,1){4.}}
\put(-2,1.151){\line(6,1){4.}}
\put(-2,1.149){\line(6,1){4.}}

\color{black}

\linethickness{0.2mm}
\put(-3,-4){\line(1,0){6}}
\put(-3, 2){\line(1,0){6}}
\put(-3,-4){\line(0,1){6}}
\put(3,-4){\line(0,1){6}}
\put(-3,1){\line(6,1){6.02}}

\multiput(-3,1)(0.25,0){24}
{\line(1,0){0.2}}

\linethickness{0.2mm}
\put(-3.8,-3){\vector(1,0){8.5}}     
\put(4.95, -3,0){\makebox(0,0){$x'$}}

\put(-0.4,  0,1){$P$}
\put(-0.8,  -2.85){$Q=(0,0)$}                
\put(0.1, -1,5){$h=|PQ|$}

\color{blue}    
\put(-1.6, 2.3){$\Gamma_{\rho}^{+}$}

\color{red}    
\put(-1.6, 1.5){$\Gamma_{\rho}^{\phi,+}$}

\color{green}    
\put(0.4, 0.3){$\Gamma_{\rho}^{\cup,+}$}
\color{black}

\put(-1.6, 1.1){${_\phi}$}
\put(2.1, -3.2){$\rho$}
\put(3.1, -3.2){$\sigma$}
\put(0.1,2.2){$_{h+2k \sigma^{1+\alpha}}$}
\put(0.1,1.2){$_{h+k \sigma^{1+\alpha}}$}
\put(0.1,-3.8){$_{-k \sigma^{1+\alpha}}$}

\end{picture}
\vspace{5.2cm}
\caption{ {\it A 2d  depiction of the cylinder      $C_{\rho}$ (blue) as well as the top surface of $\mathcal{C}_{\rho}^{\cup}$   (green) and $\mathcal{C}_{\rho}^{\phi}$(red). The bottom surface of all three domains is the same.   }}
\end{figure}
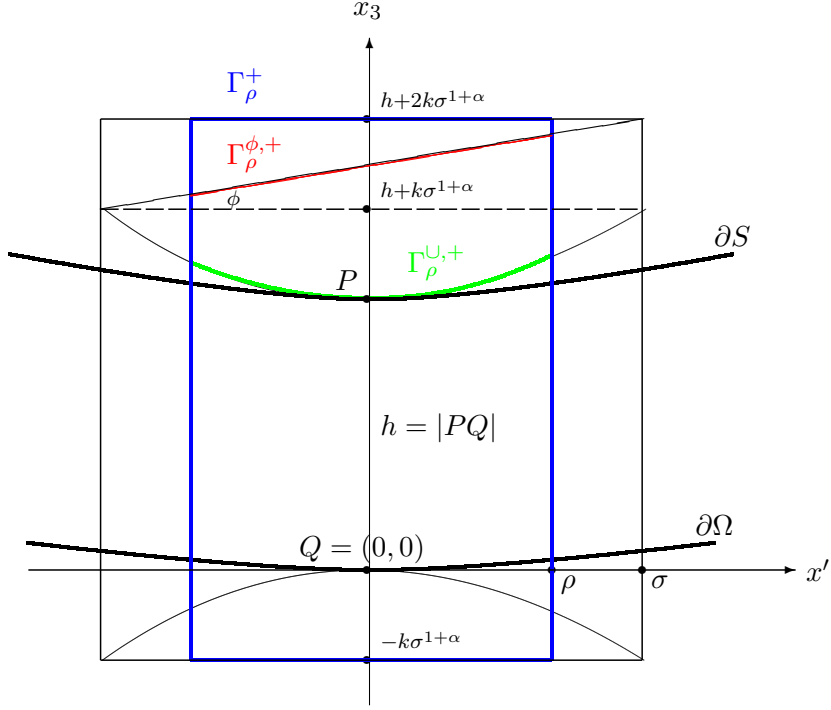

\vspace{5mm}
   
    \noindent \underline{{\bf Estimate of $\xo_{3}$.}} We will use the same argument as before, but in  a modified domain that allows for  the appearance of $\xo_{3}$. To simplify the calculations we choose our coordinate system $x_1x_2x_3$ so that $\bom = (0, \xo_2, \xo_3)$. We next
     cut the cylinder  $\mathcal{C}_{\rho, \xg }$, $\rho \leq \xs$ with the plane 
     \[
     x_3 = (x_2 - \xs)\tan \phi  + h + 2k \xs^{1+\xa},
     \]
      and we denote the resulting domain by 
   \[
\mathcal{C}_{\rho, \xg }^{\phi} := \{(r, \theta, x_3):~  (r,~ \theta) \in  D_{\rho, \xg}, ~-k \xs^{1+\xa} < x_3 < (r \sin \theta - \xs) \tan \phi +h+ 2k \xs^{1+\xa}  \}.  
 \] 
 We also denote the various parts of the boundary of $ \mathcal{C}_{\rho, \xg }^{\phi}$ by $\Gamma_{\rho, \xg}^{\phi,+}$, 
$\Gamma_{\rho,\xg}^{\phi, 0} $, 
 $\Gamma_{\rho,\xg}^{\phi, 1}$, $ \Gamma_{\rho, \xg}^{\phi, -}= \Gamma_{\rho, \xg}^{ -} $
 and $\delta \Gamma^{\phi,1}_{\rho, \xg}$ as usual.
  
 The angle  $\phi$ is so chosen that $\tan \phi = \frac{k \xs^{1+\xa}}{2 \xs} = \frac12 k\xs^{\xa}$. In particular we have 
 \[
 \mathcal{C}_{\rho, \xg }^{\cup} \subset \mathcal{C}_{\rho, \xg }^{\phi} \subset \mathcal{C}_{\rho, \xg }~~~~\mbox{ with }~~~
  \Gamma_{\rho,\xg}^{\phi, +} \subset S.
  \]
    We will apply the divergence  Theorem in $ \mathcal{C}_{\rho, \xg }^{\phi} \setminus   \mathcal{C}_{\rho/2, \xg }^{\phi} $ .
 Again, the integral over the top surface of $ \mathcal{C}_{\rho, \xg }^{\phi}$ is computed explicitly,
   \bean
 \int_{ \Gamma_{\rho, \xg}^{\phi,+}} \bu \cdot \bn \, dS  & = &
  \frac{\pi \rho^2}{2 \cos \phi} (u_{P2} \sin \phi + u_{P3} \cos \phi) +
 \frac{2}{3} \rho^3 \xo_2  \sin \xg  - \frac{2}{3} \rho^3 \xo_3  \tan \phi \sin \xg .
 \eean
For the remaining  integrals  we note that 
\[
\Gamma_{\rho,\xg}^{\phi, 0}
\subset \Gamma_{\rho,\xg}^{ 0}~~~~\mbox{and~ similarly }~~~~  \delta \Gamma^{\phi,1}_{\rho, \xg} \subset \delta \Gamma^{1}_{\rho, \xg}.
\]
 We then obtain the analogue of (\ref{a3}) or (\ref{a32})  which now reads
 \bean
 A_1 \rho^{3}\, \tan \phi  \,\sin \xg \, |\xo_3  | & \leq  &
  A_2 \rho^2 | u_{P3} | \cos \phi  +   A_2 \rho^2 | u_{P2} | \sin \phi +
  A_3  \rho^3 \sin \xg \, |\xo_2| \nonumber \\
  & & + \int_{ \Gamma_{\rho, \xg}^{0}\cup \Gamma_{\rho/2, \xg}^{0}} |\bu| dS 
 +\int_{ \delta \Gamma^{1}_{\rho, \xg} \cup \delta \Gamma^{1}_{\rho, \pi+\xg}} |\bu|\, dS,
 \eean
 with positive constants $A_i$, $i=1,2,3$, depending only on $k$, $\xa$.
  Integrating  the above estimate from  $\xg=0$ to $\xd \leq \pi$
 and then from
 $\rho=\xs/2$ 
  to $\xs$,  and employing the previous estimates, we get
\be\la{bf4}
 |\xo_3| \leq C  \xs^{-2-\xa-\frac{2}{p}} (h + 3 k \xs^{1+\xa})^{2-\frac{1}{p}}
\| \nabla \bu \|_{L^p(\cC_{\xs})} ,
\ee  
with $C=C(k,\xa)$.
   
   \vspace{5mm}
   
    \noindent \underline{\bf Completion of the proof.} Putting together (\ref{bf1}),
 (\ref{bf2})  and  (\ref{bf4})  we have shown that 
  there exists a  positive constant $C_w=C_w(k, \xa)$ depending at most  on $k$, $\xa$, s.t. 
   for any $h \in[0,H]$  and any $\xs \in (0, \xs_0/2]$  the inequalities of the statement hold true.
%
%
%
%

\finedim
 
 \vspace{3mm}
 
\noindent {\bf Remark  1.}  It is an immediate consequence of Theorem  \ref{th1} and (\ref{vic1}) that for $h \leq H$ and $\xs \in (0, \xs_0/2]$ there holds
\bea\la{cw}
 |\bu_{P}| & \leq  & C_w  \xs^{-1-\xa-\frac{2}{p}} (h + 3 k \xs^{1+\xa})^{2-\frac{1}{p}}
\| \nabla \bu \|_{L^p(\cC_{\xs})} ,  \nonumber  \\
|\bom| &  \leq  & C_w  \xs^{-2-\xa-\frac{2}{p}} (h + 3 k \xs^{1+\xa})^{2-\frac{1}{p}}
\| \nabla \bu \|_{L^p(\cC_{\xs})},
\eea
for a positive constant  $C_w=C_w(k, \xa)$. This provides an  upper estimate on the kinetique condition of the rigid body as it approaches the boundary.
 
%
%

 \vspace{3mm}
 
 \noindent {\bf Remark  2.} 
 In the special case where 
 $\xa=1$ and $p=2$ and under the assumptions of Theorem \ref{th1} we have that for a positive constant  $C_w=C_w(k)$,
\bean
| u_{P3}|& \leq & C_w \xs^{-2} (h + 3 k \xs^{3})^{\frac{3}{2}}
\| \nabla \bu \|_{L^2(\cC_{\xs})}, \\
 |\bu_{P\tau}|,  |\bom_{\tau}|   & \leq  & C_w  \xs^{-3} (h + 3 k \xs^{2})^{\frac{3}{2}}
\| \nabla \bu \|_{L^2(\cC_{\xs} )} , \\
 |\xo_{3}| &  \leq  & C_w  \xs^{-4} (h + 3 k \xs^{2})^{\frac{3}{2}}
\| \nabla \bu \|_{L^2(\cC_{\xs})}.
\eean

To prove Theorem \ref{wk}  we first present the following Lemma
\begin{lemma}\la{cyl}
We  assume that 
both $\xO$ and $S$ are uniformly $C^{1,\xa}$.
Let $\bu \in W^{1,p}_0(S, \,\xO) $ be such that ${\rm div} \bu = 0$.
There exists  positive constants $ H_{*} \leq H$ and $\xs_{*} \leq  \frac{\xs_0}{2}$ depending only on $\xa, k$, $\xs_0, H, p$ such that for
\[
h \leq H_{*} ~~~~~\mbox{and}~~~~~  \frac{\xs_0}{2} \left(\frac{h}{H} \right)^{\frac{1}{1+\xa}}  \leq \xs \leq \xs_*,
\]
there holds
\[
\int_{\cC_{\xs}} |\nabla \bu|^p dx   \leq  2   \int_{\cC_{\xs} \cap F} |\nabla \bu|^p dx .
\]

\end{lemma}

\noindent {\em Proof:} Let $h<H$ and $\xs \leq \frac{\xs_0}{2}$. For later  use  we also require require 
that
\be\la{lu}
h \xs^{-1-\xa} \leq  H \left(\frac{2}{\xs_0} \right)^{1+\xa} .
\ee
 All conditions are satisfied provided that
\be\la{inc}
h \leq H_{} ~~~~~\mbox{and}~~~~~  \frac{\xs_0}{2} \left(\frac{h}{H} \right)^{\frac{1}{1+\xa}} \leq \xs \leq \frac{\xs_0}{2}.
\ee
We next compute
\bean
\int_{\cC_{\xs}} |\nabla \bu|^p dx  & =  & \int_{\cC_{\xs} \cap F} |\nabla \bu|^p dx + \int_{\cC_{\xs} \cap S} |\nabla \bu|^p dx    \\
& = & \int_{\cC_{\xs} \cap F} |\nabla \bu|^p dx + 2^{\frac{p}{2}} |\bom|^p |\cC_{\xs} \cap S|   \\
& \leq & \int_{\cC_{\xs} \cap F} |\nabla \bu|^p   dx +  2^{\frac{p}{2}} |\bom|^p \pi \xs^{3 +\xa} (2k +k).
\eean
To continue, we  first use (\ref{cw}) and then (\ref{lu})  to get
\bean
\int_{\cC_{\xs}} |\nabla \bu|^p dx   & \leq  &  
 \int_{\cC_{\xs} \cap F} |\nabla \bu|^p dx    \\
 &  & \hspace{- 1,8 cm} +   \,  
 2^{\frac{p}{2}} \, 3  \, k \, \pi  \,
  C_w^p \, \xs^{\xa p} (h \xs^{-1-\xa} + 3k)^{2p-1} \int_{\cC_{\xs}} |\nabla \bu|^p dx  \\
&  & \hspace{- 2,4 cm} \leq  \,  \int_{\cC_{\xs} \cap F} |\nabla \bu|^p dx    \\
&  & \hspace{- 1,8 cm}  +   \, \, 2^{\frac{p}{2}} \, 3  \, k \,  \pi \,  C_w^p \,  \xs^{\xa p} \left[H \left(\frac{2}{\xs_0} \right)^{1+\xa}  \!\!\!   + 3k \right]^{2p-1}  
  \int_{\cC_{\xs}} |\nabla \bu|^p dx
\\
&  & \hspace{- 2,4 cm} = \, \int_{\cC_{\xs} \cap F} |\nabla \bu|^p dx +  \xs^{\xa p} \, B \int_{\cC_{\xs}} |\nabla \bu|^p dx ;
\eean
here $B= B(p,H, \xs_0, k, \xa)$ and is given by
\[
B:=  2^{\frac{p}{2}} \,3  \, k \,   \pi \, C_w^p  \left[H \left(\frac{2}{\xs_0} \right)^{1+\xa}  \!\!\!   + 3k \right]^{2p-1}  .
\]
We will choose $\xs$ small enough so that
\[
 \xs^{\xa p} \, B  \leq \frac12 ~~~\Leftrightarrow~~\xs \leq  (2  \, B)^{-\frac{1}{\xa p}}.
\]
For such a choice to be compatible with (\ref{inc}) we need to have
\[
\frac{\xs_0}{2} \left(\frac{h}{H} \right)^{\frac{1}{1+\xa}} 
 \leq (2  \, B)^{-\frac{1}{\xa p}}.
\]
This imposes a new smalness condition on $h$, namely that
\[
h \leq  \left(\frac{2}{\xs_0} \right)^{1+\xa} \!\! (2  \, B)^{-\frac{1+\xa}{\xa p}}~ H.
\]
Hence for 
\[
 H_{*}:= \min \{ H, ~ \left(\frac{2}{\xs_0} \right)^{1+\xa} \!\! (2  \, B)^{-\frac{1+\xa}{\xa p}}~ H\}, ~~~
 \xs_{*}:= \min \{ \frac{\xs_{0}}{2}, (2  B)^{-\frac{1}{\xa p}} \},
\]
and under the assumptions of the Lemma, we have that
\[
\int_{\cC_{\xs}} |\nabla \bu|^p dx  \leq  \int_{\cC_{\xs} \cap F} |\nabla \bu|^p dx + \frac12 \int_{\cC_{\xs}} |\nabla \bu|^p dx  ,
\]
and the result follows.

\finedim

\vspace{2mm}

\noindent {\bf Remark  3.} We note that the same result holds true in the 2D case where $\xO$, $S \subset \R^2$. The proof is essentially the same.

\vspace{3mm}
\noindent {\em Proof of Theorem \ref{wk} :} Our starting point are the estimates of Theorem \ref{th1}. Let $H_{*}$ be as in Lemma \ref{cyl}. For 
$h \leq H_{*}$ and $\xs= \xs_{h} :=  \frac{\xs_0}{2} \left(\frac{h}{H} \right)^{\frac{1}{1+\xa}}$  conditions of Lemma \ref{cyl} are fulfilled, therefore
\[
\| \nabla \bu\|_{L^p( \mathcal{C}_{\xs_h})} \leq 2^{\frac{1}{p}} \, \| \nabla \bu\|_{L^p( \mathcal{C}_{\xs_h} \cap F)} .
\]
On the other hand, to show the estimate of $|u_{P3}|$,  elementary calculations show that
\bean
\xs^{-1-\frac{2}{p}} (h + 3 k \xs^{1+\xa})^{2-\frac{1}{p}}  &  &    \\
&  & \hspace{-4cm}  = \left(\frac{\xs_0}{2  H^{\frac{1}{1+\xa}}} \right)^{1+ 2 \xa- \frac{3 + \xa}{p}} \left[H \left(\frac{2}{\xs_0} \right)^{1+\xa} + 3 k \right]^{2 -\frac{1}{p}}  
   h^{\frac{1+2 \xa}{p(1+\xa)}
  \left( p - \frac{3+\xa}{1+ 2 \xa} \right)} ,
\eean
 from which the estimate for $|u_{P3}|$ follows from the corresponding estimate of Theorem \ref{th1}. We similarly treat the other estimates of Theorem \ref{wk}.

 \finedim

\noindent {\bf Remark  4.}  In the special case where 
 $\xa=1$ and $p=2$ and under the assumptions of Theorem \ref{wk} we have that 
%
\bean
| u_{P3}| & \leq & C_w h^{\frac12} \,
\| \nabla \bu \|_{L^2(\cC_{\xs_h}\cap F)}, \\
 |\bu_{P\tau}|, |\bom_{\tau}|  & \leq  & C_w   \,
\| \nabla \bu \|_{L^2(\cC_{\xs_h}\cap F)} , \\
 |\xo_{3}| &  \leq  & C_w  h^{-\frac12}  \,
\| \nabla \bu \|_{L^2(\cC_{\xs_h}\cap F)}.
\eean

\vspace{3mm}

\noindent {\bf Remark  5.} In the 2D case  (cf Theorem 4.6 in \cite{FT}) the analogous estimates of
 Theorem \ref{wk} are
\bean  |u_{P2} | & \leq &  \, C_w  h^{\frac{1+2 \xa}{p(1+\xa)}
  \left( p - \frac{2+\xa}{1+ 2 \xa} \right)}  
\| \nabla \bu \|_{L^p(\cC_{\xs_h}\cap F)},  \\
  |u_{P1}| & \leq  & C_w  \, h^{\frac{1}{p}
  \left( p - \frac{2+\xa}{1+ \xa} \right)}  \| \nabla \bu \|_{L^p(\cC_{\xs_h}\cap F)},  \\
  |\xo|  & \leq  &  C_w  \,  h^{\frac{2\xa}{p (1+\xa)}
  \left( p - \frac{2+\xa}{ 2\xa} \right)}  \| \nabla \bu \|_{L^p(\cC_{\xs_h}\cap F)} ,
  \eean
where $u_{P2}$ and $u_{P1}$ are the vertical and tangential components of $\bu_P$ respecively and $\xo$ is the scalar angular velocity. In particular for $\xa=1$, $p=2$,
\bean 
|u_{P2} | &\leq   \, C_w  h^{\frac34}\| \nabla \bu \|_{L^2(\cC_{\xs_h}\cap F)},~~~\mbox{and}~~~~|u_{P1}| , |\xo| \leq 
\, C_w  h^{\frac14}\| \nabla \bu \|_{L^2(\cC_{\xs_h}\cap F)}.
  \eean

\setcounter{equation}{0}
\section{Estimates  for strong solutions   }

Here we will give the proof of Theorems \ref{str} and \ref{dynam}.
 Assume that both $S$ and $\xO$ are uniformy $C^3$.
For $\bx' = (x_1,x_2)$  let \[
x_3=h+g(\bx') ~~~\mbox{and} ~~~ x_3=g_{\xO}(\bx')
\]
 be the functions that describe $\partial S$   and $\partial \xO$ respectively, locally near $\bx'=0$. We then  have that
\[
|\nabla^{i} g(\bx')|\leq K |\bx'|^{2-i},  ~~ ~~i=0,1,2 ~~~\mbox{for}~~  |\bx'| < \xs_0,
\]
as well as 
\[
|\nabla^{3} g(\bx')|\leq K, ~~ \mbox{for}~~   |\bx'| < \xs_0,
\]
and similar estimates for $g_{\xO}$,
for a positive constant $K$ depending only on $S$, $\xO$. Without loss of generality we assume that 
\[
K \geq \max \{1,k\}.
\]
To prove Theorem \ref{str} let us first show 
 
 \begin{prop} \la{prop}
 We  assume that 
both $\xO$ and $S$ are uniformly $C^{3}$.
Let $\bu \in W^{1,2}_0(S, \,\xO) $ be such that ${\rm div} \bu = 0$ and in addition $\bu \in W^{2,2}(F)$ .
Then there exists a positive constant $C(K)$ depending only on $K$, such that   for any $ h \in [0, H]$ and
 $\xs \in (0, \xs_0/2]$ there holds
\bean 
|u_{P3}| & \leq  &  3 \pi^{-\frac12}  \, \xs^{-2} \,(h+ 2K \xs^2)^{\frac52}\,
 \| \nabla^2 \bu\|_{L^2( \mathcal{C}_{\xs}\cap F)} + C(K) \xs^2 |\bu_{P \tau}| \nonumber \\
 & & \,  + \,   C(K) \,  \xs \, (h  +2 K \xs^2)  |\bom_{ \tau}|.
\eean
\end{prop}

\noindent {\em Proof:}
 We  apply divergence Theorem in the domain $ \mathcal{C}_{\rho } \cap F$. Setting
 \bean
 \tilde{\Gamma}_{\rho}^{+} & = & \{|\bx'|< \rho,   ~~  x_3 = h+g(\bx')
  \} , 
 \eean
 we have
 \be\la{divs1}
 \int_{ \tilde{\Gamma}_{\rho}^{+}} \bu \cdot \bn \, dS  =  - \int_{ \Gamma_{\rho}^{0}\cap F} \bu \cdot \bn \, dS.
 \ee
   
 \noindent \underline{\bf Estimates on  $ \tilde{\Gamma}_{\rho}^{+}$. } Here $\bn =(1+| \nabla g|^2)^{\frac12} (-g_{x_1}, -g_{x_2}, 0)$. As usual we will write $\bu = \bu_P + \bom \times (\bx -\bx_P)     = \bu_P + \bom \times (x_1, x_2, g(x_1, x_2))$.
 \bean
 \int_{ \tilde{\Gamma}_{\rho}^{+}} \bu \cdot \bn \, dS  & = &  \int_{ |\bx'|< \rho}  \bu \cdot (- g_{x_1}, -g_{x_2}, 1) \,d \bx'  \\ 
 &  & \hspace{- 2,7 cm} = \int_{ |\bx'|< \rho} \bu_P \cdot \bn_P \, d \bx' +  \int_{ |\bx'|< \rho}  \bu_P \cdot (- g_{x_1}, -g_{x_2}, 0) \,d \bx' \\
 & & \hspace{- 2,5 cm}+\, \bom \cdot  \int_{ |\bx'|< \rho} (x_1, x_2, g) \times (- g_{x_1}, -g_{x_2}, 0)\, d \bx' + \bom \cdot  \int_{ |\bx'|< \rho} (x_1, x_2, g) \times (0,0, 1)\, d \bx'   \\
 &  & \hspace{- 2,7 cm} = \pi \rho^2 \bu_P \cdot \bn_P  - \bu_{P \tau} \cdot \int_{ |\bx'|< \rho} ( g_{x_1}, g_{x_2} ) dx_1 dx_2  \\
 & & \hspace{- 2,5 cm} + \, \bom \cdot  \int_{ |\bx'|< \rho} (gg_{x_2}, -gg_{x_1}, x_2g_{x_1} -x_1g_{x_2}) \, dx_1 dx_2 + \bom \cdot  \int_{ |\bx'|< \rho} (x_2, -x_1, 0) \, dx_1 dx_2  \\
  &  & \hspace{- 2,7 cm} =: \pi \rho^2 \bu_P \cdot \bn_P  -\bu_{P \tau} \cdot \bI_1 + \bom \cdot \bI_2 + \bom \cdot \bI_3 .
  \eean
It is easily seen that $\bI_3=0$. 
 To estimate the other integrals we note that
 for $|\bx'|=|(x_1,x_2)|$ small, we have  by Taylor expansion 
\bea\la{tay}
g(x_1, x_2) &  =  & \frac12 ( g_{x_1x_1}(0,0) x_1^2 +2g_{x_1x_2}(0,0) x_1 x_2 +
g_{x_2x_2}(0,0) x_2^2) + O(K|\bx'|^3) \nonumber \\
&  =  &  \frac{\rho^2}{2} \big[ g_{x_1x_1}(0,0) \cos^2 \theta +2g_{x_1x_2}(0,0) \cos \theta \sin \theta  \\
&  & + g_{x_2x_2}(0,0) \sin \theta^2 \big] + O(K\rho^3). \nonumber 
\eea
 Then,
 \bean
 \int_{ |\bx'|< \rho} g_{x_1} \, dx_1 dx_2 & = &   \int_{-\rho}^{\rho} g((\rho^2-x_2^2)^{1/2}, x_2)-g(-(\rho^2-x_2^2)^{1/2}, x_2) \, dx_2  \\
  & = &  2g_{x_1x_2}(0,0) \int_{-\rho}^{\rho}  x_2(\rho^2-x_2^2)^{1/2} \, dx_2  +  O(K \rho^4). 
 \eean
 The last integral is zero (the integrand is odd), we therefore end up with
 \[
 |\bI_1| \leq C(K)  \rho^4 .
 \]
Concerning the integral   $\bI_2$, we  note that 
\[
\int_{ |\bx'|< \rho} (  x_2g_{x_1} -x_1g_{x_2}) \, dx_1 dx_2 = \int_{ |\bx'|< \rho} (x_2, -x_1) \cdot \nabla g \, dx_1 dx_2=0.
\]
The last equality follows  integrating by parts.  
Using the Taylor  expansion of 
  $g_{x_1}$ and $g_{x_2}$, a similar argument as above shows that
  \[
  |\bI_2| \leq  C(K)  \rho^5.
  \]
Hence,
\be\la{gam}
 \int_{ \tilde{\Gamma}_{\rho}^{+}} \bu \cdot \bn \, dS  =  \pi \rho^2 \bu_P \cdot \bn_P  + R, ~~~~~|R| \leq   C(K) (|\bu_{P \tau}| + |\bom_{\tau}|\rho) \rho^4 .
 \ee

\vspace{3mm}

 \noindent \underline{\bf Estimates on  $ \Gamma_{\rho}^{0}\cap F$. } 
In what follows we will change variables by 
\be\la{chv}
\xi = \frac{x_3-g_{\xO}}{h+g - g_{\xO}}, ~~~~~~\tilde{\bu}(r, \theta, \xi) = \bu(r, \theta, x_3) .
\ee
\bean
\int_{ \Gamma_{\rho}^{0}\cap F} \bu \cdot \bn \, dS & = &  \rho \int_{0}^{ 2\pi } \int_{g_{\xO}}^{h+g} \bu \cdot \bn \, dx_3 \, d \theta =  \rho \int_{0}^{2 \pi } \int_{0}^{1} \tilde{\bu} \cdot \bn \, (h+g -g_{\xO}) \, d \xi \, d \theta  \\
& = & \rho \int_{0}^{1} \int_{0}^{2 \pi } \tilde{\bu} \cdot \bn \,(h+g -g_{\xO}) \, d \theta 
\, d \xi  \\
& = & \rho \int_{0}^{1} \Phi(\rho, \xi) \, d\xi,
\eean
with
\[
\Phi(\rho, \xi) :=  \int_{0}^{2 \pi } \tilde{\bu} \cdot \bn \,(h+g -g_{\xO}) \, d \theta .
\]
We have that $\Phi(\rho,0) = 0$, therefore by mean value Theorem there exists
$\zeta_{*} \in (0,1)$ s.t. $\Phi_{\xi}(\rho,\zeta_{*}) = \Phi(\rho,1)$. We then write
\[
\Phi(\rho, \xi) = \int_{0}^{\xi} \Phi_{\zeta}( \rho, \zeta) \, d \zeta = \int_{0}^{\xi} \int_{\zeta_{*}}^{\zeta}  \Phi_{\tau \tau}(\rho,\tau ) \, d\tau \, d \zeta  + \xi \Phi(\rho,1).
\]
As a consequence
\bean
\int_{ \Gamma_{\rho}^{0}\cap F} \bu \cdot \bn \, dS & = & \rho \int_{0}^{1}
\int_{0}^{\xi} \int_{\zeta_{*}}^{\zeta}  \Phi_{\tau \tau}(\rho, \tau) \, d\tau \, d \zeta \, d \xi + \frac12 \rho \Phi(\rho, 1) ,
\eean
and therefore
\bean
|\int_{ \Gamma_{\rho}^{0}\cap F} \bu \cdot \bn \, dS|
& \leq & \rho \int_{0}^{1}
 |\Phi_{\tau \tau}(\rho, \tau) |\, d\tau  + \frac12 \rho |\Phi(\rho,1)|.
\eean
We estimate the first term in the rhs,
\bean
\rho \int_{0}^{1} |\Phi_{\tau \tau}( \rho, \tau) |\,  d \tau & \leq  &  \rho \int_{0}^{ 2 \pi } \int_{0}^{1}  |\tilde{\bu}_{\tau \tau}| (h+g -g_{\xO})\, d \tau d \theta
 \\
 & = &   \rho  \int_{0}^{ 2\pi} \int_{g_{\xO}}^{h+g}  |\bu_{x_3 x_3} |\, (h+g -g_{\xO})^2 \, dx_3 d \theta \\
 & \leq  &   (h+ 2 K \rho^2)^2 \, \int_{0}^{ 2\pi} \int_{g_\xO}^{h+g}  | \bu_{x_3x_3} | \rho \, d x_3 d \theta .
\eean
Concerning $\Phi( \rho,1)$, for $\bn= (\cos \theta, \sin \theta,0)$ we have
\bean
\Phi(\rho, 1) & =  & \int_{0}^{2 \pi } \tilde{\bu}(\rho, \theta, 1) \cdot \bn \,(h+g -g_{\xO}) \, d \theta  \\
 & & \hspace{-2cm}  =\int_{0}^{2 \pi }\bu(\rho, \theta, h+g) \cdot \bn \,(h+g -g_{\xO}) \, d\theta  \\
  & & \hspace{-2cm} =\int_{0}^{2 \pi }\bu_P \cdot \bn \,(h+g -g_{\xO}) \, d\theta + \bom \cdot \int_{0}^{2 \pi } (\rho \cos \theta, \rho \sin \theta, g) \times  \bn
  \,(h+g -g_{\xO}) \, d\theta \\
   & & \hspace{-2cm} := I_1 + I_2.
\eean
To estimate $I_1$, we use (\ref{tay}) to obtain
\[
|I_1 | \leq C ( K) \,|\bu_{P \tau}| \rho^3,
\]
for a  positive constant $C$ depending only on $K$. Concerning $I_2$ we have that
\bean
 (x_1, x_2, g) \times  \bn &  = & (\rho \cos \theta, \rho \sin \theta, g) \times (\cos \theta, \sin \theta, 0)  \\
 & = &
  \,(-\sin \theta, \cos \theta, 0)\, g(\rho \cos \theta,\rho \sin \theta),
 \eean
from which it follows
\bean
|I_2 | & \leq &  2 \pi K \,|\bom_{ \tau}| \rho^2 (h  +2 K \rho^2) \\
& \leq &   2 \pi K \,|\bom_{ \tau}| \xs^2 (h  +2 K \xs^2) .
\eean
Putting everything together, and recalling that $\rho \leq \xs$, we have
\bea\la{phix}
|\int_{ \Gamma_{\rho}^{0}\cap F} \bu \cdot \bn \, dS|
& \leq  & (h+ 2K \xs^2)^2 \, \int_{0}^{ 2\pi} \int_{g_\xO}^{h+g}  | \bu_{x_3x_3} | \rho \, d x_3 d \theta  \nonumber \\
& & + C ( K) \,|\bu_{P \tau}| \xs^4
+  2 \pi K  \,|\bom_{ \tau}| \xs^3 (h  +2 K\xs^{2}) .
\eea

\vspace{3mm}
 \noindent \underline{\bf The estimate on $|\bu_{P3}|$ :}
   Combining (\ref{divs1}), (\ref{gam}) and (\ref{phix})  we get
\bean
\pi \rho^2 |u_{P3}| & \leq  & (h+ 2K \xs^2)^2 \, \int_{0}^{ 2\pi} \int_{g_\xO}^{h+g}  | \bu_{x_3x_3} | \rho \, d x_3 d \theta  \\
& & + \, C ( K) \,|\bu_{P \tau}| \xs^4
+  C(K)  \,|\bom_{ \tau}| \xs^3 (h  +2 K \xs^2) .
\eean
We integrate this from $\rho=0$ to $\rho=\xs$ to get after simplifying,
\bean
|u_{P3}| & \leq  &  \, \frac{3 \xs^{-3}  }{\pi} (h+ 2K \xs^2)^2 \, 
 \| \nabla^2 \bu\|_{L^1( \mathcal{C}_{\xs}\cap F)} + C(K) \xs^2 |\bu_{P \tau}| \\
 & & \, + \,   C(K) \,  \xs \, (h  +2 K \xs^2)  |\bom_{ \tau}|.
\eean
By Swchartz inequality 
\bean
 \| \nabla^2 \bu\|_{L^1( \mathcal{C}_{\xs}\cap F)} & \leq  &  |\mathcal{C}_{\xs}\cap F|^{\frac12}  \,
 \| \nabla^2 \bu\|_{L^2( \mathcal{C}_{\xs}\cap F)}  \\
  & \leq  &  \sqrt{\pi}  \, \xs (h  +2 K \xs^2)^{\frac12}\, \| \nabla^2 \bu\|_{L^2( \mathcal{C}_{\xs}\cap F)}.
 \eean
Hence
\bean 
|u_{P3}| & \leq  &  3 \pi^{-\frac12}  \, \xs^{-2} \,(h+ 2K \xs^2)^{\frac52}\,
 \| \nabla^2 \bu\|_{L^2( \mathcal{C}_{\xs}\cap F)} + C(K) \xs^2 |\bu_{P \tau}| \nonumber \\
 & & \,  + \,   C(K) \,  \xs \, (h  +2 K \xs^2)  |\bom_{ \tau}|,
\eean
and the proof is complete.
\finedim

 \vspace{3mm}
 
\noindent {\bf Remark 1 }  Examination of the proof shows that in case both $g$ and $g_{\xO}$ are radially symmetric,  the terms involving $|\bu_{P \tau}|$ and $|\bom_{ \tau}|$ are absent. Thus we recover 
the result of Proposition 8 of \cite{HT2}. However, in the  absence of symmetries these terms are present.
To proceed then with the estimation of $|u_{P3}|$ we will use the estimates  for $|\bu_{P \tau}|$ and $|\bom_{ \tau}|$ from Theorem \ref{th1} (the case of weak solutions).

\vspace{3mm}

\noindent {\bf Remark 2 } In the 2D case, under similar assumptions there exists a positive constant $C(K)$ depending only on $K$, such that   for any $ h \in [0, H]$ and
 $\xs \in (0, \xs_0/2]$ there holds
\bean 
|u_{P2}| & \leq  &  \sqrt{2}\, \xs^{-\frac32}  \,(h+ 2K \xs^2)^{\frac52}\,
 \| \nabla^2 \bu\|_{L^2( \mathcal{C}_{\xs}\cap F)} + C(K) \xs^2 |\bu_{P \tau}| \nonumber \\
 & & \,  + \,   C(K) \,  \xs^2 \, (h  +2 K \xs^2)  |\bom_{ \tau}|.
\eean

\vspace{3mm}

\noindent {\it Proof of Theorem \ref{str}}: Let  $h < H$  and  $\xs \in (0, \xs_0/2]$.
 We combine the estimate of Proposition \ref{prop} with the estimates for the weak solutions for 
$|\bu_{P \tau}|$ and $|\bom_{ \tau}|$,  with $p=2$, $\xa=1$, cf Remark 2 of previous section, namely
\[
|\bu_{P \tau}|, ~|\bom_{ \tau}|  \leq C_{w}(k) \,  \xs^{-3} \,(h+ 3k \xs^2)^{\frac32}\,
 \| \nabla \bu\|_{L^2( \mathcal{C}_{\xs})} .
 \]
Taking also into account that $K \geq k$ we obtain
\bean
|u_{P3}| & \leq &  C(k) \, \xs^{-5} \, \,(h+ 3K \xs^2)^{4}\,
 \| \nabla^2 \bu\|_{L^2( \mathcal{C}_{\xs}\cap F)} \\
 & & \hspace{-1cm} 
 +  \, C(k,K)\,  \xs^{-1} \,(h+ 3K \xs^2)^{\frac32} \,
 \| \nabla \bu\|_{L^2( \mathcal{C}_{\xs})}\\
& & \hspace{-1cm}   + \, C(k,K)\,  \xs^{-2} \,(h+ 3K \xs^2)^{\frac52} \,
 \| \nabla \bu\|_{L^2( \mathcal{C}_{\xs})} .
\eean
At this point we  further  restrict $h$ by $h \leq H_{*} \leq H$ where $H_{*}$ is defined in Lemma \ref{cyl}.
In addition we  choose $\xs=\xs_{h}:=    \left(\frac{h}{H} \right)^{\frac{1}{2}} \frac{\xs_0}{2}$.
With these choices the assumptions of Lemma \ref{cyl} are satisfied, therefore
\[
\| \nabla \bu\|_{L^2( \mathcal{C}_{\xs_h})} \leq \sqrt{2} \, \| \nabla \bu\|_{L^2( \mathcal{C}_{\xs_h} \cap F)} .
\]
On the other hand,
\[
(h+ 3K \xs^2_h) = \xs^{2}_h( h \xs^{-2}_h + 3K) 
\leq \xs^2_h  \left[H \left(\frac{2}{\xs_0} \right)^{2}  + 3K 
\right].
\]
Hence,
\[
|u_{P3}| \leq   C_{s} \left( \xs_h^3 \,
 \| \nabla^2 \bu\|_{L^2( \mathcal{C}_{\xs_h}\cap F)} +  \xs_h^2 \,
 \| \nabla \bu\|_{L^2( \mathcal{C}_{\xs_h}\cap F)}  \right) ,
 \]
with a positive constant $C_{s} = C_s (k, \xs_0, H, K)$. Recalling the definition of $\xs_h$ the result follows.

\finedim

\vspace{3mm}

\noindent {\bf Remark 3 } In the 2D case, under similar to Theorem \ref{str} assumptions, there exist positive constants $H_{*}$, $C_s$, $c_0$ depending only on $S$, $\xO$, such that, whenever $h < H_{*}$ and $\xs_{h}:= c_0 h^{\frac{1}{2}}$,
\[
|u_{P2}| \leq   C_{s} \left( h^{\frac74} \,
 \| \nabla^2 \bu\|_{L^2( \mathcal{C}_{\xs_h}\cap F)} +  h^{\frac54} \,
 \| \nabla \bu\|_{L^2( \mathcal{C}_{\xs_h}\cap F)}  \right) .
 \]

\vspace{3mm}

We next give the proof of Theorem \ref{dynam}.

\vspace{3mm}

\noindent {\it Proof of Theorem \ref{dynam}:} It is a consequence of Theorems  \ref{wk}, \ref{str} and Lemma 4.5 of \cite{FT}. We recall from \cite{FT}  that if $R$ is the radius of the inner sphere of $\xO$ and $h(t) <R$ then
 $h(t)$  is Lipschitz continuous and for almost all $t$
\[
\Big| \frac{dh(t)}{dt} \Big| = |\bu_{P} \cdot \bn_{P}| = |u_{P3}|.
\]
Therefore if $h < H_0 :=\min \{R, H_{*}\}$ the result follows.

\finedim

As a consequence of Theorem \ref{dynam}(ii) we have the following noncollision result for strong solutions.

\vspace{2mm}

\begin{coro} \la{cor}  We assume that both $\xO$ and $S_0$ are uniformly $C^{3}$  domains. Let the vector field $\bu$ be such that 
 \[
\bu \in L^{\infty}(0,T; L^{2}(\xO)) \cap L^{1}(0,T; W^{1,2}_0(S(t),\xO)),~~~ {\rm div} \bu = 0,
\]
and in addition 
\[
\int_{0}^{t} \| \nabla^2 \bu\|_{L^2( F(s))} ds < \infty, ~~~~~\forall t \in [0, T).
\]
If  $h(0)>0$  then $h(t)$ remains positive for $t \in [0,T)$. Moreover, if ~$\liminf_{t \to T} h(t)=0$, then 
\[
\lim_{t \to T}\int_{0}^{t}  h^{\frac12}(\tau) \|\nabla^2 \bu \|_{L^2( \mathcal{C}_{\xs_h(\tau)}\cap F(\tau))} \, d \tau = \infty .
\]
\end{coro}

\noindent {\it Proof}: By Theorem \ref{dynam}(ii) we have that for a.a. $t$,
\[
\left| \frac{1}{h} \frac{d h}{dt} \right| \leq C_{s} \left( h^{\frac12} \,
 \| \nabla^2 \bu\|_{L^2( \mathcal{C}_{\xs_h}\cap F)} +  
 \| \nabla \bu\|_{L^2( \mathcal{C}_{\xs_h}\cap F)}  \right) .
 \]
Integrating this we get
\[
\left|\ln \frac{h(t)}{h(0)} \right| \leq C_s \int_0^t   h^{\frac12} \|\nabla^2 \bu \|_{L^2( \mathcal{C}_{\xs}\cap F)} \, d \tau +  C_s \int_0^t \|\nabla \bu \|_{L^2( \mathcal{C}_{\xs}\cap F)} \, d \tau ,
\]
from which the result follows easily.

\finedim

\vspace{1mm}

\noindent {\bf Remark 4 } Similar results to Theorem \ref{dynam} as well as to  Corollary \ref{cor} hold  true in the 2D  case. In particular, under the assumptions of Theorem \ref{dynam} the analogue of \eqref{th19} is 
\[
\left| \frac{d h}{dt} \right|  \leq \, C_w  h^{\frac{1+2 \xa}{p(1+\xa)}
  \left( p - \frac{2+\xa}{1+ 2 \xa} \right)} 
\| \nabla \bu \|_{L^p(\cC_{\xs_h}\cap F)}, 
\]
whereas the analogue of \eqref{s4} is
\[
\left| \frac{d h}{dt} \right|  \leq \, C_{s} \left( h^{\frac74} \,
 \| \nabla^2 \bu\|_{L^2( \mathcal{C}_{\xs_h}\cap F)} +  h^{\frac54} \,
 \| \nabla \bu\|_{L^2( \mathcal{C}_{\xs_h}\cap F)}  \right) .
\]

\setcounter{equation}{0}
\section{ Example and optimality of Theorem \ref{wk}}

Let $\xO$ be the half space $\R^3_{+} = \{ (x_1, x_2, x_3): ~ x_3>0 \}$.
For  $r=|\bx'|=|(x_1, x_2)|$,
 the lower part of the  boundary of the body $S$ is given by  $x_3= h(t)+ r^{1+\xa} $ so that $h(t)$ equals the distance between $S$ and $\partial \xO$ at each time $t$. Let
\[
\bF = \frac{1}{2} (- \dot{h} x_2, ~\dot{h} x_1,  -\xo_3 r^2) ,
\]
so that the velocity of the solid is given 
by 
\[
\bu_S = \nabla \times \bF = (-\xo_3 x_2, \xo_3 x_1, \dot{h}) = (0, 0, \dot{h}) + (0,0, \xo_3) 
\times (x_1, x_2, x_3-h) .
\]
 This represents  rotation about the $x_3$ axis with angular velocity $\xo_3$ and a vertical speed $\dot{h}$. To define the velocity field everywhere in $\xO$ we set
 \bean
 \bu & =  & \nabla \times (\bF \Psi) \\
 & = &  \frac12 ( - \xo_3 x_2 (2  \Psi + r \Psi_r ) - \dot{h} x_1  \Psi_{x_3},  
 ~\xo_3 x_1 (2  \Psi + r \Psi_r ) - \dot{h} x_2  \Psi_{x_3}, ~ \dot{h}  (2  \Psi + r \Psi_r ) ) ,
\eean
so that ${\rm div} \bu = 0$.
Here 
\[
\Psi=\Psi(r, x_3)= \phi \left(\frac{x_3}{h+ r^{1+\xa}} \right)    \psi_1(r) \psi_2(x_3) ,
\]
and $\phi$, $\psi_1$, $\psi_2$ are suitable smooth cut off functions such that
\be  \nonumber
\phi(\tau)  = \left\{ \begin{array}{ll}
1,    & \tau \geq 1, \\
0, &  \tau \leq 0  \\
\end{array} ,\right.~~
 \psi_1(r)  = \left\{ \begin{array}{ll}
1,    &  r < \rho , \\
0, &  r >  2 \rho \ \\
\end{array} , \right. ~~
 \psi_2(x_3)  = \left\{ \begin{array}{ll}
1,    &  x_3 < H , \\
0, &  x_3 > 2H  \ \\ 
\end{array} . \right. 
\ee

By standard arguments we calculate for a small but fixed $\xg$.
\bean
\int_{\R^{3}_{+}} |\nabla \bu |^2 d \bx & \lesssim  &   \int_{0}^{\infty} \int_{0}^{h+  r^{1+\xa}} \dot{h}^2   r^3 \Psi_{x_3x_3}^2+ \xo_3^2 r^3 (\Psi_{x_3} +r\Psi_{rx_3})^2 \, dx_3  \,  \, dr  \\
& \lesssim  &  \int_{0}^{\xg} \int_{0}^{h+  r^{1+\xa}}\dot{h}^2  r^3 \phi_{x_3x_3}^2 
+ \xo_3^2 r^3 (\phi_{x_3} +r\phi_{rx_3})^2 \,\, dx_3 dr + O_{h}(1)  \\
&  \lesssim  &    \int_{0}^{\xg} \frac{ \dot{h}^2  r^3}{(h+  r^{1+\xa})^3} 
+ \frac{ \xo_3^2  r^3}{h+  r^{1+\xa}} \,
dr + O_{h}(1).
\eean
Hence,
\bean
\int_{\R^{3}_{+}} |\nabla \bu |^2 d \bx  \lesssim  
 \left\{ \begin{array}{ll}
\dot{h}^2 + \xo_3^2  ,   ~~~~ &  0< \xa < 1/3, \\
 \dot{h}^2  \,  |\ln h|+ \xo_3^2, & \xa = 1/3,\ \\
 \dot{h}^2  \, h ^{\frac{1- 3\xa}{1+ \xa}}+ \xo_3^2, &  1/3< \xa \leq1 \ .\\
\end{array} \right.
\eean
We similarly  have
\be\la{para}
\int_{ \mathcal{C}_{\xs_{h}}\cap F} |\nabla \bu |^2 d \bx \lesssim \dot{h}^2  \, h^{\frac{1- 3\xa}{1+ \xa}}+ \xo_3^2 \, h^{\frac{ 3-\xa}{1+ \xa}} .
\ee
\[
\int_{\R^{3}_{+}} |\bu |^2 d \bx \lesssim  \dot{h}^2 + \xo_3^2,
\]

\noindent  {\bf Optimality of the  estimates of Theorem \ref{wk}} Taking $\xo_3=0$  we have  from (\ref{para}) that 
\[
| \dot{h} |  \gtrsim    h^{\frac{ 3\xa-1}{2(1+ \xa)}}   \| \nabla \bu \|_{L^2(\cC_{\xs_h}\cap F)}, 
\]
which shows the optimality of the $u_{P3}$ estimate of Theorem \ref{wk}. Similarly for $\dot{h}=0$ we get the optimality of the $\xo_3$ estimate. Starting with the vector field 
\[
\bF = (0, \, -\frac12 \xo_2 (x_1^2 + (x_3-h)^2), \, x_2 v_1(t) ),
\]
so  that
\[
\bu_S = \nabla \times \bF = (v_1,0,0)+ (0, \xo_2,0) \times (x_1, x_2, x_3-h),
\]
 one can show the optimality of the rest of the estimates.

\vspace{3mm}

\noindent  {\bf Non zero collision speed}  In this framework, it is natural to require $\bu$ to belong  to the energy space of weak solutions of the solid--fluid interaction problem,  that is
\be\la{ener}
\bu \in L^{\infty}(0,T; L^{2}(\xO)) \cap  L^{2}(0,T; W^{1,2}_0(S(t),\xO)) .
\ee
We easily see that $\bu \in L^{\infty}(0,T; L^{2}(\xO)$ whenever 
\be\la{conom}
|\dot{h}| + |\xo_3| < C ~~~~\hbox{ for} ~~~~~t \in (0,T],  
\ee
whereas for $\xa \in (1/3, 1]$ one has  $\bu \in L^{2}(0,T; W^{1,2}_0(S(t),\xO))$ provided that
\be\la{condh}
\int_{0}^{T} \dot{h}^2 \, h ^{\frac{1- 3\xa}{1+ \xa}} \, dt < \infty ,~~~~~~~~1/3< \xa \leq1.
\ee
We take $h(t) = (T-t)^{\theta}$, $\theta \geq 1$ and $\xo(t)$ arbitrary bounded function. One can show that such a $\bu$ 
solves  the problem of  the body-- fluid interaction with a forcing term in $L^{2}(0,T; H^{-1}(\xO))$, see e.g \cite{V1,FT}.
Condition (\ref{condh}) is satisfied provided if $\theta > \frac{1 + \xa}{3 -\xa}$. Thus for $\xa<1$ we can choose $\theta=1$ and therefore we have a non zero speed of collision ( $\dot{h}=1$). On the other hand if $\xa=1$ then $\theta>1$ and we have a   zero  collision  speed, in agreement with Theorem 3.2 in \cite{V1}.  In connection with this, 
we recall that the assumption that the rigid body is smooth seems to lead to various results that are in disagreement
with experimental observation, see e.g \cite{GH1, H, HT1},  as well as \cite{LG} p. 324.
%
%
%
%
%


\begin{thebibliography}{MMM}





\bibitem {AB}
Alvarado R., Brigham D., Mazya V., Mitrea M.,  Ziade E.,
 On the regularity of domains satisfying a uniform hour-glass condition and a sharp version of the Hopf-Oleinik boundary point principle.  {\em Probl. Mat. Anal.} {\bf 57}, 3--68 (2011).  
  
%
%
%
%
%
%




\bibitem{FT} Filippas S., Tersenov A.,
 On vector fields describing the 2d motion of a rigid body in a viscous fluid and applications. 
{\em J. Math. Fluid Mech.} {\bf 23(1)} 24 p. (2021). 




%


\bibitem{GH1} G\'{e}rard--Varet D., Hillairet M.,  Computation of the drag force  on a sphere close to a wall. The Roughness issue. {\em ESAIM: M2AN},  {\bf 46}, 1201–1224,  (2012) 

\bibitem{H05} Hillairet M., Do Navier Stokes Equations Enable to Predict Contact Between Immersed Solid Particles, {\em  Advances in Mathematical Fluid Mechanics}, 109-127 (2007) 

\bibitem{H} Hillairet M.,
 Lack of collision between solid bodies in a 2D incompressible viscous flow, {\em Commun. Partial Differ. Equations}, {\bf 32(9)}, 1345--1371 (2007)
 
 \bibitem{HT1} Hillairet M., Takahashi T., Collisions in 3D Fluid Structure interaction problems, {\em  SIAM J. on Math. Analysis},  {\bf 40(6)} 2451--2477, (2009)
 
 \bibitem{HT2} Hillairet M., Takahashi T., Blow up and grazing collision in viscous fluid solid interaction systems,  {\em  Ann. I. H. Poincare  A.N.},  {\bf 27}, 291--313,  (2010)
 

%
\bibitem{HS} Hoffman K.--H., Starovoitov V. N., On a motion of a solid body in a
viscous fluid. Two dimensional case, {\em Adv. Math. Sci. Appl.}, {\bf 9}(2), 633--648, (1999).
%
\bibitem{LT} Lacave Ch.,  Takahashi T.,
Small Moving Rigid Body into a Viscous
Incompressible Fluid, {\em Arch. Rational Mech. Anal.}, {\bf 223}, 1307--1335, (2017).



\bibitem{LG} Leal, L. G.,
Advanced transport phenomena. Fluid mechanics and convective transport processes. 
Cambridge University Press, xix, 912 p. (2007).



%

\bibitem {V1} Starovoitov V. N., Behavior of Rigid Body in an Incompressible Viscous Fluid near a Boundary, {\em Int. Series of Numerical Math.},  {\bf 147}, 313--327. (2003).

\bibitem {SST} San Martin J.A., Starovoitov V. N., Tucsnak M., Global weak solutions for the two--dimensional motion of several rigid bodies in an incompressible viscous fluid, {\em Arch. Rational Mech. Anal.}, {\bf 161}, 113--147, (2002).

\bibitem {Tak}  Takahashi T.,
 Analysis of strong solutions for the equations modeling the motion of a rigid-fluid system in a bounded domain.  {\em Adv. Differ. Equ.} {\bf 8(12)}, 1499-1532, (2003). 
%
%








\end{thebibliography}
\end{document}